# Algebraic Constructions for Digraph Routing Problems

Nyumbu Chishwashwa* Vance Faber† Noah Streib‡

**Keywords:** Graph, Network Routing, Digraph, Vertex Transitive, Groupoid, Quasigroup, Spanning Factorization.
**Mathematics Subject Classification:** 05C20, 05C25, 05C38, 05C90.

**Abstract**

Efficiency of routing on a regular digraph often involves finding opitmal properties of the graph. For example, the diameter of a digraph is the maximum distance between any two vertices. We show how we can study these problems algebraically in terms of quasigroups, 1-factors, and permutation groups. Our investigation originated from the study of graphs as the Cayley graphs of groupoids with $d$ generators, a left identity, and right cancellation; that is, a right quasigroup. This enables us to provide compact algebraic definitions for some important graphs that are either given as explicit edge lists or as the Cayley coset graphs of groups larger than the graph. One such example is a single expression for the Hoffman-Singleton graph. From there, we notice that the groupoids can be represented uniquely by a set of disjoint permutations and we explore the consequences of that observation.

## Roadmap

Here is a brief roadmap to the sections and examples in this paper.

1. Groupoids and Cayley digraphs
2. Partial groupoid tables and their extensions
3. Factorizations and groupoids
4. Spanning factorizations and vertex transitive digraphs
5. Examples of Cayley digraphs of groupoids
6. Relationships between Factorization Properties
   (a) Example 6.1: Illustration of equivalences between digraph properties
   (b) Example 6.2: Kautz digraph as a Cayley digraph of a groupoid
   (c) Example 6.3: Groupoid not satisfying Axiom 1
   (d) Example 6.4: Groupoid with Cayley graph isomorphic to the Hoffman-Singleton Graph
   (e) Example 6.5: Alegre digraph as a Cayley digraph of a groupoid
7. Cyclic difference digraphs
8. Generalized cyclic difference digraphs

*University of Western Cape, Cape Town, South Africa
†Center for Computing Sciences, Bowie, Maryland
‡Center for Computing Sciences, Bowie, Maryland

9. Some properties of generalized cyclic digraphs
   (a) Example 9.3: Alegre digraph as a cyclic difference digraph
   (b) Example 9.6: A second representation of Alegre as a cyclic difference digraph
10. Covering groups of digraphs
   (a) Example 10.3: Covering group of Alegre has diameter 23 and 187,500 vertices
   (b) Example 10.4: Generalization of Alegre to GCD with 49 vertices
   (c) Example 10.5: Generating $S_3 \, wr \, S_3$ as the cover group of two permutations in $S_9$
   (d) Example 10.10: Three digraphs with degree 2, diameter 2, and 6 vertices
   (e) Example 10.11: Three digraphs with degree 2, diameter 3, and 12 vertices

# 1 Groupoids and Cayley Digraphs

Here, we define a *groupoid* as a finite set $\Gamma$ equipped with a binary (not necessarily associative) operation $*$, referred to as the *product* (often suppressed in notation). When we index the elements of the groupoid, the matrix of products is called the *table* of the groupoid.

A *word* $\omega$ in the elements of $\Gamma$ is a finite sequence of elements, and the *value* of $\omega$ is the element formed by taking products from left to right in $\Gamma$. A subset $S \subset \Gamma$ is called a *generating set* for $\Gamma$ (denoted $\Gamma = \langle S \rangle$) if each element in $\Gamma$ is the value of some word in $S$.

Given a groupoid $\Gamma = \langle S \rangle$, we create an associated digraph $G(S)$ with directed edges $(u, u * s)$ for each $u \in \Gamma$ and $s \in S$. By analogy to group theory, we call this the *Cayley digraph* of the groupoid.

# 2 Partial Groupoid Tables and Their Extensions

We want our Cayley digraphs to be regular and loopless. What does this imply for the groupoid? In a Cayley digraph of a groupoid generated by a set of size $d$, the $d$ columns corresponding to the generators encode the $d$ edges in the graph emanating from each vertex $u$, as the edges are exactly $(u, u * s)$. We call these $d$ columns a *partial groupoid*.

If, for some element $x$ in the groupoid, we have $x * s = x$, then the Cayley digraph will have a loop. Since we are interested in loopless digraphs, we forbid our groupoids from containing such elements; that is, they are *loopless* groupoids: $x \notin xS$ for all $x \in \Gamma$. If two partial groupoids have isomorphic Cayley graphs, we say they are *equivalent*. Any groupoid with the same columns is called an *extension* of the partial groupoid.

Given any connected, regular, and loopless digraph $G$, we will demonstrate how to construct a partial groupoid $\Gamma$ with generating set $S$ such that $G$ is the Cayley digraph of $\Gamma = \langle S \rangle$, and $\Gamma$ has the following properties:

(P1) The empty word, denoted $e$, is a left identity.

(P2) $x \notin x * S$ for all $x \in \Gamma$ (loopless).

(P3) $u * s = u * t$ for $u \in \Gamma$ and $s, t \in S$ only if $s = t$ (left cancellation).

We will continue to refer to this structure as a groupoid, even with these additional properties.To achieve this, we first need to develop some necessary tools.

## 3 Factorizations and groupoids

A *factor* of a regular digraph $G$ is a function (written as applying from the right instead of the usual left) on the vertices which maps each vertex to a directed neighbor. We also think of a factor as a set of edges of $G$. A factor is a *1-factor* if it is one-to-one. A set of factors $F = \{F_1, F_2, \ldots, F_d\}$ of a regular digraph $G$ of degree $d$ is a *factorization* if the factors are pairwise disjoint. Let $F_1$, $F_2$, …, $F_d$ be the factors in a factoring of $G$. A *word* $\omega$ in the elements of $F$ is a finite sequence of elements, and the *value* of $\omega$ is the function formed by taking compositions from left to right in $F$. If $v$ is a vertex and $\omega$ is a word, then $v\omega$ denotes the directed path (and its endpoint) in $G$ starting at $v$ and proceeding along the unique edge corresponding to each consecutive factor represented in the word $\omega$. We say a set of words $W$ is *tree-like* if it is prefix-closed: whenever $\omega \in W$, every initial subword of $\omega$ also lies in $W$. By Petersen's theorem (see, for example, [LST92]), every regular digraph has a 1-factorization. For completeness, we provide a proof.

*Theorem* (Petersen). Every loopless digraph $G$ where the in-degree and out-degree of every vertex is $d$ has an edge-disjoint decomposition into $d$ 1-factors.

*Proof.* Form an auxiliary graph $B$ with two new vertices $u'$ and $u''$ for each vertex $u$. The edges of $B$ are the pairs $(u', v'')$ where $(u, v)$ is a directed edge in $G$. The undirected graph $B$ is bipartite and regular with degree $d$, and so by Hall's Marriage Theorem, it can be decomposed into $d$ 1-factors. Each of these 1-factors corresponds to a directed 1-factor in $G$. □

Now we can prove a fundamental theorem (see [Dör74]).

**Theorem 3.1** (Factorizations $\iff$ Cayley graphs of groupoids)**.** *Fix an integer $d \geq 1$. There is a natural bijection, up to isomorphism and permutation of labels, between*

- *isomorphism classes of pairs $(G, F)$ where $G$ is a connected, loopless, $d$–regular digraph and $F = \{F_1, \ldots, F_d\}$ is a factorization of $G$ into $d$ disjoint 1–factors (written on the right), and*
- *isomorphism classes of partial groupoids $(\Gamma, S)$ with $S = \{s_1, \ldots, s_d\}$ satisfying (P1) left identity, (P2) no loops, and (P3) left cancellation, where $G$ is the Cayley digraph of $(\Gamma, S)$ with edges $u \to u * s_i$.*

*More precisely:*

(A) ***From** $(G, F)$ **to** $(\Gamma, S)$. Given $(G, F)$ and a choice of root $r \in V(G)$, there exists a tree-like set $W$ of words in the alphabet $F$ labeling the unique $r$-rooted breadth–first spanning tree. Let $\Gamma$ be the set of these words and set $S = \{s_i\}$ with $s_i$ the singleton word $F_i$. Define the product $\omega * s_i$ by concatenation; when the edge labeled $F_i$ leaves the vertex $\omega$ in $G$, $\omega * s_i$ is the word labeling its head. Then $(\Gamma, S)$ satisfies (P1)–(P3), and $G$ is isomorphic to the Cayley digraph of $(\Gamma, S)$.*

(B) ***From** $(\Gamma, S)$ **to** $(G, F)$. Given $(\Gamma, S)$ with (P1)–(P3), form the Cayley digraph $G$ on vertex set $\Gamma$ with edges $u \to u * s_i$, and set*

$$F_i = \{(u,\, u * s_i) \mid u \in \Gamma\} \quad (1 \leq i \leq d).$$

*Then each $F_i$ is a factor, the family $F = \{F_1, \ldots, F_d\}$ is a factorization, and $G$ is connected, loopless, and $d$–regular.*

*These two constructions are mutually inverse up to isomorphism: different choices of the root $r$ in (A) or permutations of the labels $\{1, \ldots, d\}$ induce isomorphic groupoids (with relabeled generators) and the same isomorphism class of Cayley digraphs, and conversely relabeling $S$ corresponds to permuting the factors in $F$.*

*Proof.* We first check direction (A). Given $(G, F)$ with $G$ connected, loopless, and $d$–regular and $F = \{F_1, \ldots, F_d\}$ a factorization, fix a root $r \in V(G)$. The breadth–first tree $T$ rooted at $r$ uniquely labels each vertex by a prefix–closed word $\omega$ in the factors, with edges corresponding to appending singletons $F_i$. Declaring $S = \{s_i\}$ with $s_i$ the generator corresponding to $F_i$, and defining $\omega * s_i$ to be the word labeling the neighbor along $F_i$, produces a partial groupoid $(\Gamma, S)$. This structure satisfies (P1)–(P3): the empty word is a left identity, looplessness ensures no $u * s_i = u$, and left cancellation holds because each $F_i$ is a 1–factor. By construction, the Cayley digraph of $(\Gamma, S)$ coincides with $G$.

For direction (B), suppose $(\Gamma, S)$ is a partial groupoid with (P1)–(P3). Form $G(S)$ with edges $u \to u * s_i$ whenever defined, and let $F_i = \{(u, u * s_i) \mid u \in \Gamma\}$. Then each $F_i$ is a 1–factor: every vertex has exactly one outgoing edge in $F_i$ by definition, and at most one incoming edge by left cancellation (P3). Thus $G(S)$ is $d$–regular and loopless, and strong connectivity follows because every word in $S$ labels a directed path.

Finally, the two constructions are inverse up to isomorphism. Different root choices or relabelings of $S$ change only the labeling of generators, not the underlying isomorphism class. □

**Definition 3.2.** Using the construction in Theorem 3.1, we can form a *canonical* extension of the partial groupoid $\Gamma$ with properties (P1)–(P3) to a groupoid. By construction, every element $\omega \in \Gamma$ is represented uniquely by a finite word in the generators $S = \{s_1, \ldots, s_d\}$. Given $\mu \in \Gamma$ and a word $\omega = s_{i_1} s_{i_2} \cdots s_{i_k}$, we define the product $\mu * \omega$ recursively:

$$\mu * s_{i_1} := \text{the element given by the partial product in } \Gamma,$$

and for $k > 1$,

$$\mu * (s_{i_1} s_{i_2} \cdots s_{i_k}) := (\mu * s_{i_1}) * (s_{i_2} \cdots s_{i_k}).$$

This defines a total operation of right multiplication by any word $\omega$ in the generators, extending the partial groupoid product to all of $\Gamma$. The result is a groupoid in which multiplication by generators agrees with the original partial structure.

*Remark* 3.3. It should be noted that (P1) is independent of the other two properties by examining a simple example (see the particular assignment for the 6-vertex Kautz graph in Example 6.3). In many cases, we want our groupoids to have an additional property called *right cancellation*:

(P4) $u * s = v * s$ for $u, v \in \Gamma$ and $s \in S$ only if $u = v$.

The properties (P3) and (P4) were defined on a partial groupoid so they most naturally correspond to the Cayley graph which is solely determined by the partial groupoid. These properties do not always extend to the canonical extension. We differentiate explicitly by defining the following two additional properties:

(P5) $u * v = u * w$ for $u, v, w \in \Gamma$ only if $v = w$ (left quasigroup).

(P6) $u * w = v * w$ for $u, v, w \in \Gamma$ only if $u = v$ (right quasigroup).

## 4 Spanning factorizations and vertex transitive digraphs

In this section, we review what we know about spanning factorizations [DF14].

**Definition 4.1.** Let $G$ be a regular digraph with disjoint factors $F_1, \ldots, F_d$. For each vertex $v$, let $W(v)$ be a set of words in the factors, containing the empty word $\omega_0 = e$ and the singleton words $\omega_i = F_i$ for $1 \leq i \leq d$, such that the set of endpoints $vW(v)$ covers all vertices of $G$. The collection $\{W(v) \mid v \in V(G)\}$ forms an *all-to-all routing* on $G$. If there exists a single set of words $W$ with $W(v) = W$ for every vertex $v$, then we say that the factorization (or 1-factorization) together with $W$ forms a *spanning factorization* (respectively, a *spanning* 1-*factorization*) of $G$.

Note that in this definition, we require the words $\omega_i$ for $1 \leq i \leq d$ to be the factors $F_i$. This requirement is absent in [DF14].

Note also that given a vertex $v$, a word $\omega$ in $W(v)$ and a factor $F_i$ in $\omega$, then $F_i$ represents a unique edge in $G$.

**Definition 4.2.** A spanning factorization is *neighborhood preserving* if for every $\alpha, \beta \in W$ and $F_i$, we have

$$(\alpha\beta)F_i = \alpha(\beta F_k)$$

for some $k$ (note $k$ can depend on $\alpha, \beta$, and $F_i$).

**Definition 4.3** (Vertex transitive)**.** A digraph $G$ is *vertex transitive* if for any two vertices $u$ and $v$ there is an automorphism of $G$ which maps $u$ to $v$.

We will also need to use the Cayley coset representation of a vertex transitive digraph.

**Definition 4.4** (Cayley coset graph)**.** Let $\Gamma$ be a finite group, $\mathrm{H}$ a subgroup, and $S$ a subset. Suppose the following hold:

1. $S \cap \mathrm{H} = \emptyset$ and $\Gamma$ is generated by $S \cup \mathrm{H}$,
2. $\mathrm{H}S\mathrm{H} \subseteq S\mathrm{H}$,
3. $S$ is a subset of distinct coset representatives of $\mathrm{H}$ in $\Gamma$.

Then we can form the *Cayley coset digraph* $G = (\Gamma, S, \mathrm{H})$ with the cosets $\{g\mathrm{H} : g \in \Gamma\}$ as vertices and the set of pairs $(g\mathrm{H}, gs\mathrm{H})$ with $s \in S$ as edges. When $\mathrm{H}$ is the identity subgroup, the graph is a *Cayley digraph*.

The classic proof of Sabidussi [Sab64] shows that a digraph is vertex transitive if and only if it is a Cayley coset digraph. An important aspect of the proof shows that one can construct a Cayley coset digraph from a vertex transitive digraph by using the automorphism group as the group $\Gamma$ required in the definition and the subgroup of automorphisms that fix a vertex as the required subgroup $\mathrm{H}$. The generators $S$ correspond to automorphisms that map a vertex to a neighbor. In [CFS24], we investigated whether or not a regular digraph $D$ has a spanning factorization. In particular, we focused on vertex transitive digraphs of degree 2. We were able to show that using the original definition in [DF14], they all do. Here, we extend that investigation to Cayley graphs of groupoids.

## 5 Relationships between Factorization Properties

**Definition 5.1** (Notation)**.** We assign acronyms to the properties of a connected digraph as follows:

**(VT)** vertex transitive,

**(SF)** has a spanning factorization,

**(S1F)** has a spanning 1-factorization,

**(TL)** is SF and tree-like,

**(NP)** is SF and neighborhood preserving.

In [CFS24], we showed the following theorems.

**Theorem 5.2.** *(VT $\Rightarrow$ SF, TL, NP) A connected vertex transitive digraph $G = (\Gamma, S, H)$ has a tree-like, spanning factorization which is neighborhood preserving.*

**Theorem 5.3.** *(SF, NP $\Rightarrow$ VT) A digraph $G$ with a neighborhood preserving spanning factorization*

$$W = \{\omega_0 = \emptyset,\ \omega_1 = F_1, \ldots, \omega_d = F_d,\ \omega_{d+1}, \ldots, \omega_n\}$$

*is vertex transitive.*

**Theorem 5.4.** *(SF, NP $\Leftrightarrow$ VT) A digraph $G$ is vertex transitive if and only if it has a neighborhood preserving spanning factorization.*

**Example 5.5.** (S1F $\not\Rightarrow$ VT) There exists a digraph that has a spanning 1-factorization but is not vertex transitive. The graph has vertices $0, \ldots, 7$ modulo 8 and degree 2. The edges are the cycle $(i, i+1)$, the four edges $(i, i+2)$ for odd $i$, and the four edges $(i, i+4)$ for even $i$. The 1-factor $F_1$ is the cycle; the remaining edges form $F_2$. The spanning 1-factorization is then

$$W = \{\emptyset, F_1, F_2, F_1F_2, F_2F_1, F_1^2, F_2F_1^2, F_1F_2F_1, F_2F_1F_2\}.$$

**Example 5.6.** (VT $\not\Rightarrow$ S1F+TL) There exists a vertex transitive digraph of degree 2 which has no tree-like spanning 1-factorization. This digraph is $G = (\Gamma, S, H)$ with $\Gamma = A_5$, $S = \{(0,1,2,3,4), (0,3,4,2,1)\}$ and $H = (0,2)(1,3)$. See [Fab25] for details.

### Groupoid interpretations of the properties

Each of the properties in Definition 5.1 for a graph $G$, corresponds to properties in some groupoid whose Cayley graph is $G$.

- **(SF)** This is the content of (P5).
- **(S1F)** This is (P5) together with (P4).
- **(TL)** This means the words in $W$ are prefix-closed.
- **(NP)** Neighborhood preserving corresponds to quasi-associativity of the groupoid: for any elements $u, v$ and generator $s_i$ there exists a generator $s_j$ such that $(u * v) * s_i = u * (v * s_j)$.
- **(VT)** This is equivalent to **SF** together with **NP** (see Mwambene [Mwa06]).

### Open Questions

The digraph in Example 1 is not tree-like. We cannot rule out the implications:

$$\text{VT} \Rightarrow \text{S1F}, \qquad \text{SF, TL} \Rightarrow \text{VT}, \qquad \text{S1F, TL} \Rightarrow \text{VT}.$$

*Remark* 5.7. a) The words in $F$ and $S$ are applied in opposite order, $u\omega(F) = \omega(S)(u)$. See Example 6.1 below.
b) The earliest use of groupoids to create vertex transitive graphs seems to be [Gau96], which adds an extra condition that the vertices correspond to a set of automorphisms with the property that for each pair of vertices $u$ and $v$, there is exactly one automorphism in the set mapping $u$ to $v$.
c) The words generated by $F$ are permutations on the vertices, and so are the words generated by the automorphisms in $S$. Even though we showed that these are associated with each other when the graph is vertex transitive, they are not the same group. See Example 6.1. Below, we will come back to the group generated by the 1-factors.

## 6 Examples of Cayley Digraphs of Groupoids

We provide examples to illustrate how groupoids can describe digraphs.

**Example 6.1.** This example illustrates the equivalences of Theorem 5.4. The group is a semidirect product group on the two generators $\alpha$ and $\theta$ with

$$\alpha^4 = 1 = \theta^5, \quad \alpha^{-1}\theta\alpha = \theta^2.$$

A permutation representation is given by

$$\theta = (0,1,2,3,4)(5,6,7,8,9), \quad \alpha = (0,5)(1,8,4,7)(2,6,3,9).$$

The subgroup is $H = \{1, \alpha^2\}$. The vertices of the graph are the $H$-cosets of

$$[e, \theta, \theta^2, \theta^3, \theta^4, \alpha, \theta\alpha, \theta^2\alpha, \theta^3\alpha, \theta^4\alpha] = [0,1,2,3,4,5,6,7,8,9].$$

The edges and 1-factors corresponding to $\theta$ and $\alpha$ are computed using $\alpha^{-1}\theta\alpha = \theta^2$:

$$F_1 = \{(0,1,2,3,4),(5,8,6,9,7)\}, \quad F_2 = \{(0,5),(1,6),(2,7),(3,8),(3,9)\}.$$

The group generated by the automorphisms $\alpha$ and $\theta$ is not the same as the group generated by the permutations $F_1$ and $F_2$. The latter group is isomorphic to the semidirect product $\mathbb{Z}_2(\mathbb{Z}_5 \times \mathbb{Z}_5)$.
The words $\omega_i$ forming the tree of paths starting from 0 are given by:

$$\omega_0 = \emptyset \to \omega_5 = F_2,$$

$$\begin{aligned}
\omega_1 = F_1 &\to \omega_6 = F_1F_2,\\
\omega_2 = F_1^2 &\to \omega_7 = F_1^2F_2,\\
\omega_3 = F_1^3 &\to \omega_8 = F_1^3F_2,\\
\omega_4 = F_1^4 &\to \omega_9 = F_1^4F_2.
\end{aligned}$$

It is easy to check that $i = 0\omega_i(F_1, F_2) = \omega_i(\theta, \alpha)(0)$. We can now complete the $\omega_1 = F_1$ and $\omega_2 = F_2$ columns in the groupoid multiplication table, corresponding to the graph edges.

| | $\omega_1$ | $\omega_5$ |
|---|---|---|
| $\omega_0$ | $\omega_1$ | $\omega_5$ |
| $\omega_1$ | $\omega_2$ | $\omega_6$ |
| $\omega_2$ | $\omega_3$ | $\omega_7$ |
| $\omega_3$ | $\omega_4$ | $\omega_8$ |
| $\omega_4$ | $\omega_0$ | $\omega_9$ |
| $\omega_5$ | $\omega_8$ | $\omega_0$ |
| $\omega_6$ | $\omega_9$ | $\omega_1$ |
| $\omega_7$ | $\omega_5$ | $\omega_2$ |
| $\omega_8$ | $\omega_6$ | $\omega_3$ |
| $\omega_9$ | $\omega_7$ | $\omega_4$ |

**Example 6.2.** This example is a groupoid on $\mathbb{Z}_2 \times \mathbb{Z}_3$. Multiplication is defined by the table

| | 00 | 01 | 02 | 10 | 11 | 12 |
|---|---|---|---|---|---|---|
| 00 | 00 | 01 | 02 | 10 | 11 | 12 |
| 01 | 01 | 02 | 10 | 12 | 00 | 01 |
| 02 | 02 | 10 | 11 | 01 | 02 | 10 |
| 10 | 10 | 11 | 12 | 00 | 01 | 02 |
| 11 | 11 | 12 | 00 | 02 | 10 | 11 |
| 12 | 12 | 00 | 01 | 11 | 12 | 00 |

Note that the columns are permutations, but the rows are not. The generators are $t = (1, 0)$ and $s = (0, 1)$; these columns are fixed-point-free. The generator $s$ produces a 6-cycle in the Cayley graph. The set $H = \{e, t\}$ is a subgroup, and $Hs$ and $Hs^2$ are disjoint cosets. In the Cayley graph, $H$ forms a 2-cycle, and $Hs \cup Hs^2$ forms a 4-cycle. This graph is the Kautz graph $G(2, 3)$.

**Example 6.3.** The following groupoid on $\mathbb{Z}_2 \times \mathbb{Z}_3$ satisfies the second and third axioms but not the first. It is finitely generated by $t = (1, 0)$ and $s = (0, 1)$, but the empty word cannot be assigned a consistent meaning. The issue is that it acts as an identity on the right but not on the left, demonstrating that axiom 1 is independent of the other axioms.

| | 00 | 01 | 02 | 10 | 11 | 12 |
|---|---|---|---|---|---|---|
| 00 | 00 | 01 | 02 | 11 | 12 | 10 |
| 01 | 01 | 02 | 00 | 10 | 11 | 12 |
| 02 | 02 | 00 | 01 | 12 | 10 | 11 |
| 10 | 10 | 11 | 12 | 01 | 02 | 00 |
| 11 | 11 | 12 | 10 | 00 | 01 | 02 |
| 12 | 12 | 10 | 11 | 02 | 00 | 01 |

Even though this is not a groupoid in the sense we are using, it still has a Cayley graph, which is the same directed graph as in Example 2.

**Example 6.4.** We define a groupoid on $\mathbb{Z}_2 \times \mathbb{Z}_p \times \mathbb{Z}_p$ by

$$(a, b, c) * (x, y, z) = (a + x, b - bx + y, c + (-1)^a by + 2^a z)$$

for $p$ a prime. For $p = 5$,

$$S = \{(0, 0, 1), (0, 0, 4), (1, 0, 0), (1, 1, 0), (1, 2, 0), (1, 3, 0), (1, 4, 0)\}$$

produces the Hoffman-Singleton graph [HS60], an undirected graph of degree 7 and diameter 2.

**Example 6.5.** The Alegre graph [FYA84], the largest known degree 2 diameter 4 digraph, has 25 vertices. It can be represented as the Cayley graph $G$ of a groupoid on $\mathbb{Z}_5 \times \mathbb{Z}_5$. The elements are $t^i s^j$, represented as natural numbers from 0 to 24 in lexicographic order. The factor determined by $t$ consists of a 5-cycle starting at 0, another 5-cycle starting at 3, and a 15-cycle:

$$(0, 5, 10, 15, 20)(3, 23, 18, 13, 8)(1, 17, 24, 21, 12, 19, 16, 7, 14, 11, 2, 9, 6, 22, 4).$$

## 7 Cyclic difference digraphs

*Note.* Up to this point, a word in the 1-factors has been parsed from left to right with vertices on the left. From now on, we will want to think of a 1-factor as a permutation of the vertices and therefore parse words from right to left and apply them as functions with domain elements on the right. So, given a 1-factor $F$, it is a permutation with an edge out from $v$ being $(v, F(v))$. We give a general construction of a family of digraphs of degree 2 which have some symmetry but are not necessarily vertex transitive. We will call the members of this family *cyclic difference digraphs*. (These are similar to the metacyclic graphs in [AP82] but those are designed to be vertex transitive.) We start with vertices $V = \{k : 0 \le k < n\}$ and assume that $n = ab$ with both $a$ and $b$ not equal to 1. Each cyclic difference digraph is the disjoint union of two special 1-factors, $Z$ and $Y$. It is convenient to think of a 1-factor as both a set of edges and a permutation on $n$ so an edge in the 1-factor $F$ is $(k, F(k))$. The permutation $Z$ is a single cycle, $Z(k) = k + 1 \pmod{n}$. We call an interval of the form $[ia, (i + 1)a)$ the $i$ *segment*. To construct the complementary 1-factor $Y$, we specify every element of $V$ by giving its segment and the position in the segment:

that is, $k = ia + j$ with $0 \le j < a$ and $0 \le i < b$. Let $T = \{t_0, t_1, \ldots, t_{a-1}\}$ be a sequence of (not necessarily distinct) elements of $[0, b)$ and let $\pi$ be a permutation on $[0, a)$. The edges in $Y$ are defined for each $k = ia + j$ by

$$(ia + j, (i + t_j)a + \pi(j))$$

where the second entry is taken modulo $n$. To ensure that $Y$ and $Z$ form complementary 1-factors, we need to impose additional conditions on the sequence $T$ and the permutation $\pi$. We will find these conditions in the next section where we generalize this definition.

# 8 Generalized cyclic difference digraphs

In this section, we give another construction of degree 2 digraphs. Such a graph has two disjoint 1-factors. We can call these 1-factors $F_1$ and $F_2$ and we overload the symbol $F_i$ to denote either 1) all the edges in the 1-factor, 2) the edge $vw$, or 3) a function which inputs a vertex $v$ and outputs the edge $vw$ in $F_i$. The intended meaning should be clear from context. Note that in the context of 3), the two 1-factors are disjoint derangements on the set of vertices. As shown above in Theorem 3.1, each degree 2 digraph corresponds to a groupoid with two columns that are disjoint derangements. We create a digraph on the Cartesian product $V = \mathbb{Z}_a \times \mathbb{Z}_b$ where $n = ab$. We can write each element in $V$ as $k = (j, i)$ where $j \in \mathbb{Z}_a$ and $i \in \mathbb{Z}_b$. Our construction generalizes the vertex transitive metacirculant digraphs defined in [AP82].

**Lemma 8.1.** *Let $(X_j \mid j \in \mathbb{Z}_a)$ be a collection of permutations on $\mathbb{Z}_b$ and let $\sigma$ be a permutation on $\mathbb{Z}_a$. If for each $j$, either $X_j$ is a derangement or $j$ is not a fixed point of $\sigma$, the function $X(j, i) = (\sigma(j), X_j(i))$ is a derangement on $V = \mathbb{Z}_a \times \mathbb{Z}_b$.*

*Proof.* Assume that $(\sigma(j), X_j(i)) = (\sigma(j'), X_{j'}(i'))$. Then in particular, $j = j'$ and so $X_j(i) = X_j(i')$. This in turn implies that $i = i'$, so $X$ is a permutation. Similarly, if $(\sigma(j), X_j(i)) = (j, i)$ for some $i$ and $j$, then $\sigma(j) = j$ and $X_j(i) = i$, so $j$ is a fixed point of $\sigma$ and $i$ is a fixed point of $X_j$, contradicting the hypothesis. □

We call the function $X(j, i) = (\sigma(j), X_j(i))$ defined in Lemma 8.1 a *semi-direct permutation* on $V = \mathbb{Z}_a \times \mathbb{Z}_b$.
*Note.* Given two semi-direct permutations on $V = \mathbb{Z}_a \times \mathbb{Z}_b$, $A(j, i) = (\alpha(j), A_j(i))$ and $B(j, i) = (\beta(j), B_j(i))$, where $\alpha$ and $\beta$ are permutations on $\mathbb{Z}_a$ and the $A_j$ and $B_j$ are permutations on $\mathbb{Z}_b$, their composition is the permutation $AB$ given by $AB(j, i) = (\alpha\beta(j), A_{\beta^{-1}(j)}B_j(i))$. This is the standard product formula in the wreath product of $S_a$ and $S_b$, which is a semidirect product of $S_a \wr S_b$. Many permutations in $S_a \wr S_b$ are not derangements and therefore not eligible to be 1-factors of a graph. We will return to this wreath product later.

**Lemma 8.2.** *Given the semi-direct derangement $Z(i, j) = (\sigma(j), Z_j(i))$ and the semi-direct permutation $T(i, j) = (\theta(j), T_j(i))$, then*

1. *$Z$ and $Y = ZT$ are disjoint if and only if $T$ is a derangement,*
2. *$Y = ZT$ is a derangement if and only if the digraph formed by $Z$ and $T$ has no dicycle.*

*Proof.* First, since $Z$ is a derangement, it is a 1-factor. Then $Z(k) = Y(k) = ZT(k)$ if and only if $T(k) = k$. Second, $ZT(k) = k$ if and only if $(k, T(k))$ and $(T(k), k)$ are both edges in the union of $Z$ and $T$. □

**Definition 8.3.** We call the digraph constructed in this way a *generalized cyclic difference digraph* or GCD.

**Theorem 8.4.** *A cyclic difference digraph is a GCD.*

*Proof.* We let $Z = (\varsigma, (Z_j))$ be the derangement $Z(k) = k + 1 \pmod{ab}$ on $V$, $k = ai + j$. This is realized by $\varsigma(j) = j + 1 \pmod a$ and $Z_j(i) = i$ except $Z_{a-1}(i) = i + 1 \pmod b$. Let $T = (\theta, (T_j))$ with $\theta(j) = \pi(j) - 1 \pmod a$ and $T_j(i) = i + t_j \bmod (b)$ unless $j = \theta(a-1)$, when $T_{\theta(a-1)}(i) = i - 1 + t_{\theta(a-1)} \bmod (b)$. We then calculate

$$Y(j,i) = ZT(j,i) = (\varsigma\theta(j), Z_{\theta^{-1}(j)}T_j(i)) = (\pi(j), Z_{\theta^{-1}(j)}T_j(i)).$$

We have $Z_{\theta^{-1}(j)}T_j(i) = Z_{\theta^{-1}(j)}(i + t_j) = i + t_j$ unless $\theta^{-1}(j) = a - 1$ when

$$Z_{a-1}T_{\theta(a-1)}(i) = T_{\theta(a-1)}(i) + 1 = (i - 1 + t_{\theta(a-1)}) + 1 = i + t_{\theta(a-1)}.$$

In any case, $Y(j,i) = (\pi(j), i + t_j)$, which proves the theorem. □

**Theorem 8.5.** *A line digraph of degree 2 digraph is a GCD.*

*Proof.* We let the 1-factors of the graph $G$ be $F_0$ and $F_1$. The GCD will have vertices $\mathbb{Z}_2 \times \mathbb{Z}_n$, which we write as $(j,i)$. We let $Z_0(i) = F_0(i)$, $Z_1(i) = F_1(i)$, and $\varsigma$ be the identity. Clearly, $Z = (\varsigma, Z_j)$ is a semi-direct derangement by Lemma 8.1. We let $T(j,i) = (\theta(j), i)$ where $\theta$ is the transposition $(0,1)$, which is also a semi-direct derangement by Lemma 8.1. We claim that the GCD with factors $Z$ and $Y = ZT$ is isomorphic to the line graph of $G$. Suppose $i$ is a vertex in $G$. The edges from $i$ are $(i, F_0(i))$ and $(i, F_1(i))$, which are the vertices in the line graph. We can assign these vertices the indices $(0,i)$ and $(1,i)$, respectively. There are four types of edges in the line graph, and each one is an edge in $Z$ or $Y$:

$$((i, F_0(i)), (F_0(i), F_0^2(i))) = ((0,i), (0, F_0(i))) = Z(0,i),$$

$$((i, F_1(i)), (F_1(i), F_1^2(i))) = ((1,i), (1, F_1(i))) = Z(1,i),$$

$$((i, F_0(i)), (F_0(i), F_1F_0(i))) = ((0,i), (1, F_0(i))) = ZT(0,i),$$

$$((i, F_1(i)), (F_1(i), F_0F_1(i))) = ((1,i), (1, F_1(i))) = ZT(1,i).$$

□

## 9 Some properties of generalized cyclic difference digraphs

Let $M$ be a generalized cyclic difference digraph generated by the semi-direct derangements $Z$ and $Y = ZT$. Graph automorphisms of $M$ are elements of $S_n$ which preserve edges. Suppose $\alpha$ is an automorphism and $v$ is a vertex. Then, since the edges out from $v$ are $(v, Z(v))$ and $(v, Y(v))$, an automorphism satisfies either:
A1) $\alpha Z(v) = Z\alpha(v)$ and $\alpha Y(v) = Y\alpha(v)$
or
A2) $\alpha Z(v) = Y\alpha(v)$ and $\alpha Y(v) = Z\alpha(v)$.

**Theorem 9.1.** *A permutation $\alpha$ is an automorphism of the GCD $M$ with factors $Z$ and $Y = ZT$ if and only if for every vertex $v$, either*
*B1) $\alpha(v) = Z^{-1}\alpha Z(v) = T^{-1}Z^{-1}\alpha ZT(v)$*
*or*
*B2) $\alpha(v) = Z^{-1}\alpha ZT(v) = T^{-1}Z^{-1}\alpha Z(v)$.*

*Proof.* In case A1) $\alpha Z(v) = Z\alpha(v)$ implies $Z^{-1}\alpha Z(v) = \alpha(v)$ and $T^{-1}Z^{-1}\alpha ZT(v) = Y^{-1}\alpha Y(v) = \alpha(v)$. In case A2) $\alpha Z(v) = Y\alpha(v)$ yields $Z^{-1}\alpha Z(v) = T\alpha(v)$ and $\alpha Y(v) = Z\alpha(v)$ yields $Z^{-1}\alpha ZT(v) = \alpha(v)$. Put these together to get $T^{-1}Z^{-1}\alpha Z(v) = Z^{-1}\alpha ZT(v) = \alpha(v)$. Conversely, if B1) holds, then A1) is clear. If B2) holds, $\alpha(v) = Z^{-1}\alpha ZT(v)$ yields $Z\alpha(v) = \alpha Y(v)$, while $\alpha(v) = T^{-1}Z^{-1}\alpha Z(v)$ yields $Y\alpha(v) = \alpha Z(v)$. □

**Lemma 9.2.** *Consider a cyclic difference graph with $Y(j,i) = ZT(j,i) = (\pi(j), i+t_j)$. The cycle of the permutation $Y(j,i) = ZT(j,i) = (\pi(j), i+t_j)$ containing the element $(j,i)$ has length $\alpha c$, where $c$ is the length of the cycle of $\pi$ containing $j$, and $\alpha > 0$ is the smallest integer such that*

$$\alpha(t_j + t_{\pi(j)} + t_{\pi^2(j)} + \ldots + t_{\pi^{c-1}(j)}) = 0 \pmod b.$$

*Proof.* Start at $(j,i)$ and apply $Y$ repeatedly. We get a sequence of elements

$$(\pi^k(j), i + t_j + t_{\pi(j)} + t_{\pi^2(j)} + \ldots + t_{\pi^{k-1}(j)})$$

and if this sequence repeats $(j,i)$, we must have

$$j = \pi^k(j)$$

and

$$t_j + t_{\pi(j)} + t_{\pi^2(j)} + \ldots + t_{\pi^{k-1}(j)} = 0.$$

This can only happen when $k = \alpha c$, so

$$\alpha(t_j + t_{\pi(j)} + t_{\pi^2(j)} + \ldots + t_{\pi^{c-1}(j)}) = 0$$

and the lemma follows. □

**Example 9.3.** We can generate the Alegre graph if we choose the following values for the parameters:

$$n = 25,$$

$$a = 5,$$

$$b = 5,$$

$$\pi = (0, 2, 4),$$

$$(t_0, t_1, t_2, t_3, t_4) = (1, 4, 4, 1, 4).$$

The resulting 1-factor $Y$ is

$$(0, 7, 4, 20, 2, 24, 15, 22, 19, 10, 17, 14, 5, 12, 9)(1, 21, 16, 11, 6)(3, 8, 13, 18, 23).$$

**Lemma 9.4.** *Let $(j,i)$ be any vertex of a cyclic difference graph $G$ with $Y(j,i) = ZT(j,i) = (\pi(j), i+t_j)$. The map $\tau(j,i) = (j, i+1)$ is an automorphism of $G$.*

*Proof.* Let $k = ia + j$, so $\tau Z(k) = \tau(ia + j + 1) = ((i+1)a + j + 1)$ and

$$Z\tau(k) = Z((i+1)a + j) = (i+1)a + j + 1.$$

In addition, $\tau Y(k) = \tau(\pi(j), i+t_j) = (\pi(j), i+t_j+1)$, and $Y\tau(k) = Y(i+1, j) = (i+1+t_j, \pi(j))$. Therefore, $\tau$ is an automorphism by A1. □

**Lemma 9.5.** *The map $\mu(j) = j + 1 \pmod a$ creates an isomorphism between the cyclic difference graph with $Y(j,i) = ZT(j,i) = (\pi(j), i+t_j)$ and the cyclic difference graph with $Y'(j,i) = ZT'(j,i) = (\pi'(j), i+t'_j)$, where $\pi'(j) = \pi(j-1) + 1 \pmod a$, $t'_0 = t_{a-1} - 1$, $t'_{\pi^{-1}(a-1)+1} = t_{\pi^{-1}(a-1)} + 1$, and $t'_j = t_{j-1}$ otherwise.*

*Proof.* This isomorphism simply renames the vertex $k \in \mathbb{Z}_n$ by $k+1$ and adjusts the sequence $(t_j)$ accordingly. □

**Example 9.6.** If we apply this isomorphism twice to the generators in Example 9.3, we get a new set of parameters for an isomorphic graph:

$$V = \mathbb{Z}_{25},$$

$$n = 25,$$

$$a = 5,$$

$$b = 5,$$

$$\pi = (4, 1, 2),$$

$$(t_0, t_1, t_2, t_3, t_4) = (4, 3, 1, 1, 0).$$

The $Y$-factor corresponds to the Alegre graph as given in Example 9.3.

*Voltage digraphs.* We note a similarity between this construction and that of digraphs which are lifts of a base voltage digraph. For example, in [[DFM+19]; Figure 2], the weights in the base digraph are similar in spirit to the vector $t$ of offsets in Example 9.3 and Example 9.6. This similarity suggests a connection, though we have not yet fully identified the relationship. See also [LŠ08], where semidirect products of cyclic groups are used as voltage groups to find large digraphs with small diameters.
*Calculating the diameter.* To calculate the diameter, we only need to calculate the distance from each vertex $(j, 0)$ with $0 \leq j \leq a-1$, as shown in Lemma 9.4.

## 10 Covering groups

One of our (so far unrealized) goals is to use the algebraic machinery we constructed here to find digraphs that exceed the size of the corresponding iterated line graphs of the Alegre graph. Because we have shown in Theorem 8.4 that all these known examples are proper (neither $a$ nor $b$ is 1) generalized cyclic difference digraphs, we will restrict our search to this type of graph. In this section, we give one more tool that might help in the search.
Our starting point is Lemma 8.2, which defines a GCD. Let
$(Z_j \mid j \in \mathbb{Z}_a)$ and $(T_j \mid j \in \mathbb{Z}_a)$ be two collections of permutations on $\mathbb{Z}_b$, and let $\varsigma$ and $\nu$ be permutations on $\mathbb{Z}_a$ with the property that $Z(j,i) = (\varsigma(j), Z_j(i))$ and $T(j,i) = (\theta(j), T_j(i))$ are semi-direct derangements and $Z$ and $Y = ZT$ are disjoint 1-factors of a degree 2 digraph $M$ with $n$ vertices $V = \mathbb{Z}_b \times \mathbb{Z}_a$. More explicitly, $Y(j,i) = (\upsilon(j), Y_j(i))$ with $Y_j(i) = Z_{\theta^{-1}(j)} T_j(i)$ and $\upsilon(j) = \varsigma\theta(j)$. Paths in $M$ starting at a vertex $u$ are products of the derangements $Y$ and $Z$ applied to $u$. These products generate a subgroup $\Gamma = \langle Z, Y \rangle$ of permutations in $S_n$. To understand this group better, we consider some facts about its elements, which are permutations but clearly may not be derangements.

**Lemma 10.1.** *The semidirect permutations on $V = \mathbb{Z}_a \times \mathbb{Z}_b$ form a group $U_{ab}$ isomorphic to the semidirect product (wreath product) $S_a N$, where the normal subgroup $N = S_b^a$. This is often written as $S_b \wr S_a$. Given the permutations $A(j,i) = (\alpha(j), A_j(i))$ and $B(j,i) = (\beta(j), B_j(i))$, their composition is the semidirect permutation $(\alpha\beta(j), A_{\beta^{-1}(j)} B_j(i))$. The inverse of $A$ is $A^{-1}(j,i) = (\alpha^{-1}(j), (A_{\alpha(j)})^{-1}(i))$.*

*Proof.* We noted this product formula before in Lemma 8.2. The formula for the inverse and the fact that $N$ is normal follow as a standard exercise. □

*Notation.* We can write the semidirect permutations $(e, (A_j))$ and $(\alpha, (e_j))$ in $U_{ab}$ as $(A_j)$ and $\alpha$, respectively. This allows us to remove one set of parentheses so that $(\alpha, (A_j))$ becomes a group product $\alpha(A_j)$ in $S_a N$. Since $\alpha(A_j)\beta(B_j) = \alpha\beta(\beta^{-1}(A_j)\beta)(B_j)$, $\beta^{-1}(A_j)\beta$ is an automorphism on $N$ because it is normal and an inner automorphism on $U_{ab}$. These elements are permutations, so we can calculate that $\beta^{-1}(A_j)\beta = (A_{\beta^{-1}(j)})$, explaining where the $\beta^{-1}$ comes from.

**Definition 10.2** (Covering group)**.** Given two disjoint semidirect derangements $A$ and $B$ on $\mathbb{Z}_a \times \mathbb{Z}_b$, we call the group $\Gamma = \langle A, B \rangle$ the *covering group* in $S_{ab}$ of the digraph $G$ formed by their union. As noted above, paths in $G$ starting at a vertex $u$ are products of the derangements $A$ and $B$ applied to $u$. Thus, the covering group contains all the paths between vertices in $G$. Since all covering groups are subgroups of $U_{ab}$, we can call it the *universal covering group*.

*Note.* Given a permutation $Q \in U_{ab}$, we can find unique permutations $q \in S_a$ and $Q_j \in S_b$ such that $Q(j, i) = (q(j), Q_j(i))$ by solving $j' = q(j)$ and $i' = Q_j(i)$. The assumption that $Q \in U_{ab}$ means that $q$ must be a unique permutation independent of $i$. Additionally, if we fix $j$, the mapping of $Q$ from $i$ to $i'$ determines a unique permutation $Q_j(i)$.

**Example 10.3.** We can find the covering group $\Gamma$ of the Alegre digraph with the generators used in Example 9.3. We have the derangements

$$\rho = (0, 1, 2, 3, 4, 5, 6, 7, 8, 9, 10, 11, 12, 13, 14, 15, 16, 17, 18, 19, 20, 21, 22, 23, 24)$$

and

$$\sigma = (0, 7, 4, 20, 2, 24, 15, 22, 19, 10, 17, 14, 5, 12, 9)(1, 21, 16, 11, 6)(3, 8, 13, 18, 23).$$

Let $C_i = (i, 5 + i, 10 + i, 15 + i, 20 + i)$, $U_i = (5i + 2, 5i - 2)$, $V_i = (5i, 5i + 6)$, $a_{3i} = 20i$, $a_{3i+1} = 20i + 7$, $a_{3i+2} = 20i + 4$, $b_{5i} = 5i + 7$, $b_{5i+1} = 5i + 21$, $b_{5i+2} = 5i + 24$,

$$b_{5i+3} = 5i + 8, \quad b_{5i+4} = 5i + 20,$$

$$T = (a_0, a_1, a_2, \cdots, a_{3i}, a_{3i+1}, a_{3i+2}, \cdots, a_{12}, a_{13}, a_{14}),$$

$$\theta = (b_0, b_1, b_2, b_3, b_4, \cdots, b_{5i}, b_{5i+1}, b_{5i+2}, b_{5i+3}, b_{5i+4}, \cdots, b_{20}, b_{21}, b_{22}, b_{23}, b_{24})$$

with all values modulo 25 and $0 \le i < 5$. Let $\pi = (0, 2, 4)$ on $0 \le i < 5$. Then

1. $\rho C_i \rho^{-1} = C_{i+1}$;
2. $\sigma = T C_1^4 C_3$;
3. $\sigma^3 = (C_0 C_2 C_4)^4 C_1^2 C_3^3$;
4. $\rho^5 = C_0 C_1 C_2 C_3 C_4$;
5. $\rho^{-1}\sigma = (0, 6)(1, 20)(2, 23)(3, 7)(5, 11)(8, 12)(10, 16)(13, 17)(15, 21)(18, 22)(4, 19, 9, 24, 14) = U_0 U_1 U_2 U_3 U_4 V_0 V_1 V_2 V_3 V_4 C_4^3$;
6. $(\rho^{-1}\sigma)^2 = C_4$;
7. $\sigma\rho\sigma^{-1} = (7, 21, 24, 8, 20, 12, 1, 4, 13, 0, 17, 9, 18, 5, 22, 11, 14, 23, 10, 22, 16, 19, 3, 15) = \theta$;
8. $\theta^5 = C_0 C_1 C_2 C_3 C_4 = \rho^5$;
9. $\sigma C_i \sigma^{-1} = C_{\pi(i)}$.

There is a subgroup $Q = C_0 \otimes C_1 \otimes C_2 \otimes C_3 \otimes C_4$ in $\Gamma$. Note that from (1) $\rho C_i \rho^{-1} = C_{i+1}$. From (9) $\sigma C_i \sigma^{-1} = C_{\pi(i)}$. Thus, the action of $\Gamma$ on $Q$ is the action of the group generated by $\pi = (0, 2, 4)$ and $(0, 1, 2, 3, 4)$. This group is the alternating group $A_5$, making $\Gamma$ isomorphic to the semidirect product of $A_5$ with the elementary group $Q$. A computer search shows the diameter of the Cayley graph with generators $\rho$ and $\sigma$ is 23, with only 11 elements at this maximum distance. This group is quite possibly the largest known vertex-transitive digraph with degree 2 and diameter 23. The size of this group is $187,500 = (1.6954)^{23}$. It's challenging to determine the largest known values for such large diameters, as many constructions do not produce infinite sequences. Relevant references for vertex-transitive graphs of degree 2 include [CF], [CFG94], [DH94], [CG82], [CF95], [LŠ08],[Fab25], and [DF04].

**Example 10.4.** We can generalize this to other primes. For $p = 7$, let $C_i = (i, 7+i, 14+i, 21+i, 28+i, 35+i, 42+i)$, $U_i = (7i+2, 7i-2)$, $V_i = (7i, 7i-6)$, $W_i = (7i+4, 7i-4)$. Define $\rho$ as the cyclic permutation on $\mathbb{Z}_{49}$ and

$$\rho^{-1}\sigma = U_0U_1U_2U_3U_4U_5U_6V_0V_1V_2V_3V_4V_5V_6W_0W_1W_2W_3W_4W_5W_6C_6^4.$$

Using computational methods, we found that the group $\Gamma$ is the semidirect product $S_7\mathbb{Z}_7^7$. The diameter of the digraph is 7 compared to the diameter of the line graph of the Alegre graph, which has a diameter of 5.

**Example 10.5.** Let $A(j, i) = (\alpha(j), A_j(i))$ and $B(j, i) = (\beta(j), B_j(i))$ with $\alpha = (0, 1, 2)$, $\beta = (0, 1)$, $A_2 = (0, 1)$, $B_1 = (0, 1, 2)$, and the rest are identity permutations. Then, as permutations in $S_9$,
$A = (0, 3, 1, 4, 2, 5)$ and $B = (0, 3, 6, 1, 4, 7)(2, 5, 8)$. Using a computer, we find that these permutations generate a group $G$ with diameter 14, isomorphic to $S_3$ wr $S_3$. This group can also be generated by the two disjoint derangements $(0, 3)(1, 4)(2, 5)(6, 7, 8)$ and $(0, 7, 1, 6)(2, 8)(3, 4, 5)$.

**Theorem 10.6.** *For all integers $a, b \geq 2$, the wreath product*

$$U_{ab} \;=\; (S_b)^a \rtimes S_a \;=\; S_b \wr S_a \;\leq S_{ab}$$

*is generated by two edge–disjoint derangements in its natural action on ab points.*

*Proof.* We identify the $ab$ points with pairs $(i, j)$ where $1 \leq i \leq a$ (the "block" index) and $1 \leq j \leq b$ (the position inside the block). For $f \in S_a$ and $x \in S_b$, write the permutation

$$(i, j) \longmapsto (\, f(i),\, x(j)\, )$$

as $(x; f)$. In particular, $(x; f)$ acts by $f$ on the block index and by $x$ inside each block.

**Step 1: Choice of the top permutations $f_1, f_2 \in S_a$.** We choose $f_1, f_2$ so that $\langle f_1, f_2\rangle = S_a$ and, moreover,

$$\forall i \in \{1, \ldots, a\} \quad f_1(i) \neq f_2(i). \qquad (*)$$

Such pairs exist for every $a \geq 2$; for example:

$$\begin{array}{ll} a = 2: & f_1 = (1\,2),\ f_2 = \text{id}; \\ a = 3: & f_1 = (1\,2),\ f_2 = (2\,3); \\ a \geq 4 \text{ even}: & f_1 = (1\,2)(3\,4)\cdots(a-1\ a),\ f_2 = (1\,2\,\ldots\,a); \\ a \geq 5 \text{ odd}: & f_1 = (1\,2)(3\,4),\ f_2 = (2\,3\,\ldots\,a). \end{array}$$

**Step 2: Choice of the base permutations $x, y \in S_b$.** Choose $x, y \in S_b$ so that

$$\langle x, y\rangle = S_b. \qquad (\dagger)$$

For instance, for $b \geq 3$ take $x = (1\,2\,\ldots\,b)$ and $y = (1\,2)$; for $b = 2$ take $x = y = (1\,2)$.

**Step 3: Define the two permutations.** Set

$$\alpha = (x; f_1), \qquad \beta = (y; f_2).$$

We claim that $\alpha$ and $\beta$ are (i) derangements on the $ab$ points, (ii) edge–disjoint, and (iii) generate $S_b \wr S_a$.

*(i) $\alpha$ and $\beta$ are derangements.* A point $(i, j)$ is fixed by $(x; f)$ iff $f(i) = i$ and $x(j) = j$. By our choices, no $(i, j)$ is fixed by either $\alpha$ or $\beta$.

*(ii) $\alpha$ and $\beta$ are edge–disjoint.* For any $(i, j)$ we have

$$\alpha(i, j) = (\, f_1(i),\, x(j)\,), \qquad \beta(i, j) = (\, f_2(i),\, y(j)\,).$$

By $(*)$, $f_1(i) \neq f_2(i)$ for every $i$, so $\alpha(i, j) \neq \beta(i, j)$.

*(iii) $\langle\alpha, \beta\rangle = S_b \wr S_a$.* Let $H = \langle\alpha, \beta\rangle$. Projection to the top coordinate gives $\pi(H) = \langle f_1, f_2\rangle = S_a$. Conjugating $\beta$ places $y$ in any block, and similarly with $\alpha$ places $x$ in any block. Since $x, y$ generate $S_b$, we get the whole base group $(S_b)^a \leq H$. Thus $H = (S_b)^a \rtimes S_a = S_b \wr S_a$.

□

*Remark* 10.7. East and Mitchell [EM21] showed that $S_b \wr S_a$ is always generated by two elements, though without the additional requirement that the generators be derangements. Theorem 10.6 strengthens their result by exhibiting two *edge–disjoint derangements* that generate the wreath product, so that their union defines a 2–regular digraph with covering group the full wreath product. The case $a = 3$ and $b = 3$ in the proof was found by a computer search [Hol].

*Remark* 10.8. (One conclusion from Lemma 10.1 is that $a!(b!)^a$ divides $(ab)!$ but that is probably not relevant.) We can extend these discussions to graphs generated by three or more semidirect derangements by just imposing the condition that any pair of permutations form a GCD. We have focused on the case of degree 2 because it is the simplest case to consider. Also note that we might have worked with arbitrary derangements (1-factors) of a digraph instead of forcing them to be semidirect. In this case, the universal covering group would just be $S_{ab}$. We chose to focus on semidirect derangements because all the digraphs of record have this property.

*Remark* 10.9. Hamiltonian cycles. It seems to be an open conjecture whether every digraph of degree 2 and diameter $D$ at least as dense as the corresponding Kautz graph ($n \geq 3(2^{D-1})$) has a Hamiltonian cycle. For this reason, all of our computer searches have started with one Hamiltonian cycle. The state of the art for Hamiltonian cycles in digraphs is discussed in [KO12]. It is shown in [KO12, Theorem 30] that certain expander graphs must have Hamiltonian cycles, but we don't know if our dense graphs qualify as expander graphs of this type.

**Example 10.10.** The Kautz digraph of diameter 2 and degree 2 has 6 vertices. There are two other digraphs with the same parameters.

1. The Kautz graph of degree 2 is the line graph of the complete digraph on 3 vertices, so it is a non-trivial GCD.

2. A second digraph is the GCD given by

$$Z = (e, ((0, 1, 2), (0, 2, 1)))$$

$$Y = ((0, 1), (e, (0, 1, 2))).$$

   Written as permutations in $S_6$, these generators are

$$\rho = (0, 4, 2, 3, 1, 5)$$

$$\sigma = (0, 2, 1)(4, 5, 3),$$

   which has no cycles of length 2. This digraph is mentioned in [FYA84] where it is called $G_2^2$.

3. A third digraph of diameter 2 on 6 vertices has 1-factors

$$\rho = (0, 1, 2, 3, 4, 5)$$

$$\sigma = (0, 2, 5, 3, 1, 4).$$

This digraph has no cycles of length 2, nor does it have a 1-factor that is the union of two directed triangles, so it has a unique decomposition into 1-factors. This means that its covering group $\Gamma$ is unique, and we compute it to have order 120 and diameter 10. This graph is not a non-trivial GCD, because if it were, then its covering group would have an order that either divides $2!(3!)^2 = 72$ or $3!(2!)^3 = 48$. Examining the elements of this group shows that $\Gamma$ is one of six conjugate subgroups of $S_6$ which act transitively on the base set. This is deemed "exotic" because it only happens for a subgroup isomorphic to $S_{n-1}$ in $S_n$ when $n = 6$.

**Example 10.11.** In [MŠ01], it is shown that the largest digraph with degree 2 and diameter 3 has 12 vertices. A computer search shows that the only such digraphs are the line digraphs of the digraphs in Example 10.10. We discuss them one by one.

1. The line graph of the Kautz graph with 6 vertices is the Kautz graph with 12 vertices. It has a companion cycle $Y = (0, 6)(1, 4, 9, 8, 5)(2, 11, 7, 10, 3)$. The automorphism group is $S_3$.

2. The line graph of $G_2^2$ has a companion cycle $Y = (0, 4, 11, 9, 7, 3, 1, 6, 10, 5, 2, 8)$. The automorphism group is $\mathbb{Z}_3$.

3. The line graph of the third example above has a companion cycle $Y = (0, 10, 4, 8, 3, 11, 6, 2, 9, 1, 7, 5)$. The automorphism group is $\mathbb{Z}_4$.